\input amstex\documentstyle{amsppt}  
\pagewidth{12.5cm}\pageheight{19cm}\magnification\magstep1
\topmatter
\title Families and Springer's correspondence\endtitle
\author G. Lusztig\endauthor
\address{Department of Mathematics, M.I.T., Cambridge, MA 02139}\endaddress
\thanks{Supported in part by National Science Foundation grant DMS-0758262.}\endthanks
\endtopmatter   
\document
\define\uuE{\un{\uE}}
\define\bcp{\bar{\cp}}

\define\Irr{\text{\rm Irr}}

\define\dg{\dot g}

\define\uE{\un E}

\define\hb{\hat b}
\define\ha{\hat a}

\define\si{\sim}

\define\sqc{\sqcup}

\define\hA{\hat A}

\define\hG{\hat G}

\define\bA{\bar A}

\define\lb{\linebreak}

\define\op{\oplus}

\define\part{\partial}
\define\em{\emptyset}

\define\n{\notin}

\define\m{\mapsto}
\define\do{\dots}

\define\lra{\leftrightarrow}

\define\sub{\subset}    
\define\bxt{\boxtimes}
\define\T{\times}
\define\ti{\tilde}
\define\nl{\newline}
\redefine\i{^{-1}}

\define\un{\underline}

\define\ot{\otimes}
\define\bbq{\bar{\QQ}_l}

\define\sg{\text{\rm sgn}}

\define\a{\alpha}

\define\g{\gamma}

\define\e{\epsilon}

\redefine\o{\omega}
\define\p{\pi}
\define\ph{\phi}

\define\r{\rho}
\define\s{\sigma}

\redefine\l{\lambda}

\redefine\G{\Gamma}
\redefine\D{\Delta}

\define\kk{\bold k}

\define\NN{\bold N}

\define\QQ{\bold Q}

\define\SS{\bold S}

\define\ZZ{\bold Z}
\define\XX{\bold X}

\define\cb{\Cal B}
\define\cc{\Cal C}

\define\ce{\Cal E}
\define\cf{\Cal F}
\define\cg{\Cal G}
\define\ch{\Cal H}
\define\ci{\Cal I}
\define\cj{\Cal J}
\define\ck{\Cal K}
\define\cl{\Cal L}

\define\cp{\Cal P}

\define\cs{\Cal S}

\define\cv{\Cal V}

\define\cx{\Cal X}

\define\fa{\frak a}

\define\fA{\frak A}
\define\fB{\frak B}

\define\fT{\frak T}

\define\tA{\ti A}

\define\oca{\overset\circ\to a}
\define\ocb{\overset\circ\to b}
\define\dcup{\dot{\cup}}
\define\AL{AL}
\define\BM{BM}
\define\UE{L1}
\define\AC{L2}
\define\ORA{L3}
\define\ICC{L4}
\define\CSII{L5}
\define\CS{L6}
\define\INT{L7}
\define\SUP{L8}
\define\SH{S1}
\define\SHH{S2}
\define\SPA{Sp1}
\define\SPAA{Sp2}
\define\SP{Spr}

\head Introduction\endhead
\subhead 0.1\endsubhead
Let $G$ be a connected reductive algebraic group over an algebraically closed field $\kk$ of characteristic $p$.
Let $W$ be the Weyl group of $G$; let $\Irr W$ be a set of representatives for the isomorphism classes of 
irreducible representations of $W$ over $\bbq$, an algebraic closure of the field of $l$-adic numbers ($l$ is a 
fixed prime number $\ne p$). 

Now $\Irr W$ is partitioned into subsets called {\it families} as in \cite{\UE, Sec.9}, \cite{\ORA, 4.2}. 
Moreover to each family $\cf$ in $\Irr W$, a certain set $\XX_\cf$, a pairing $\{,\}:\XX_\cf\T\XX_\cf@>>>\bbq$, 
and an imbedding $\cf@>>>\XX_\cf$ was canonically attached in \cite{\UE},\cite{\ORA, Ch.4}. (The set $\XX_\cf$
with the pairing $\{,\}$, which can be viewed as a nonabelian analogue of a symplectic vector space, plays a key
role in the classification of unipotent representations of a finite Chevalley group \cite{\ORA} and in that of 
unipotent character sheaves on $G$). In \cite{\UE},\cite{\ORA} it is shown that $\XX_\cf=M(\cg_\cf)$ where 
$\cg_\cf$ is a certain finite group associated to $\cf$ and, for any finite group $\G$, $M(\G)$ is the set of all
pairs $(g,\r)$ where $g$ is an element of $\G$ defined up to conjugacy and $\r$ is an irreducible representation 
over $\bbq$ (up to isomorphism) of the centralizer of $g$ in $\G$; moreover $\{,\}$ is given by the "nonabelian 
Fourier transform matrix" of \cite{\UE, Sec.4} for $\cg_\cf$.

In the remainder of this paper we assume that $p$ is not a bad prime for $G$. In this case a 
uniform definition of the group $\cg_\cf$ was proposed in \cite{\ORA, 13.1} in terms of special unipotent classes
in $G$ and the Springer correspondence, but the fact that this leads to a group isomorphic to $\cg_\cf$ as 
defined in \cite{\ORA, Ch.4} was stated in \cite{\ORA, (13.1.3)} without proof. One of the aims of this paper is 
to supply the missing proof. 

To state the results of this paper we need some definitions. For $E\in\Irr W$ let $a_E\in\NN,b_E\in\NN$ be as in 
\cite{\ORA, 4.1}. As noted in \cite{\AC}, for $E\in\Irr W$ we have 

(a) $a_E\le b_E$;
\nl
we say that $E$ is {\it special} if $a_E=b_E$.

For $g\in G$ let $Z_G(g)$ or $Z(g)$ be 
the centralizer of $g$ in $G$ and let $A_G(g)$ or $A(g)$ be the group of connected components of $Z(g)$.
Let $C$ be a unipotent conjugacy class in $G$ and let $u\in C$. Let $\cb_u$ be the variety of Borel subgroups of 
$G$ that contain $u$; this is a nonempty variety of dimension, say, $e_C$. The conjugation action of $Z(u)$ on 
$\cb_u$ induces an action of $A(u)$ on $\SS_u:=H^{2e_C}(\cb_u,\bbq)$. Now $W$ acts on $\SS_u$ by Springer's 
representation \cite{\SP}; however here we
adopt the definition of the $W$-action on $\SS_u$ given in \cite{\ICC} which differs from Springer's original 
definition by tensoring by sign. The $W$-action on $\SS_u$ commutes with the $A(u)$-action. Hence we have 
canonically $\SS_u=\op_{E\in\Irr W}E\ot\cv_E$ (as $W\T A(u)$-modules) where $\cv_E$ are finite dimensional 
$\bbq$-vector spaces with $A(u)$-action. Let $\Irr_CW=\{E\in\Irr W;\cv_E\ne0\}$; this set does not depend on the 
choice of $u$ in $C$. By \cite{\SP}, the sets $\Irr_CW$ (for $C$ variable) form a partition of $\Irr W$; also, if
$E\in\Irr_CW$ then $\cv_E$ is an irreducible $A(u)$-module and, if $E\ne E'$ in $\Irr_CW$, then the 
$A(u)$-modules $\cv_E,\cv_{E'}$ are not isomorphic. By \cite{\BM} we have 

(b) $e_C\le b_E$ for any $E\in\Irr_CW$
\nl
and the equality $b_E=e_C$ holds for exactly one $E\in\Irr_CW$ which we denote by $E_C$ (for this $E$, $\cv_E$ is 
the unit representation of $A(u)$).

Following \cite{\ORA, (13.1.1)} we say that $C$ is {\it special} if $E_C$ is special. (This concept was 
introduced in \cite{\AC, Sec.9} although the word "special" was not used there.) From (b) we see that $C$ is 
special if and only if $a_{E_C}=e_C$. 

Now assume that $C$ is special. We denote by $\cf\sub\Irr W$ the family that contains $E_C$. (Note that $C\m\cf$ 
is a bijection from the set of special unipotent classes in $G$ to the set of families in $\Irr W$.) We set
$\Irr^*_CW=\{E\in\Irr_CW;E\in\cf\}$ and
$$\ck(u)=\{a\in A(u);a\text{ acts trivially on $\cv_E$ for any }E\in\Irr^*_CW\}.$$
This is a normal subgroup of $A(u)$. We set $\bA(u)=A(u)/\ck(u)$, a quotient group of $A(u)$. Now, for any 
$E\in\Irr^*_CW$, $\cv_E$ is naturally an (irreducible) $\bA(u)$-module. Another definition of $\bA_u$ is given in 
\cite{\ORA, (13.1.1)}. In that definition $\Irr^*_CW$ is replaced by $\{E\in\Irr_CW;a_E=e_C\}$ and 
$\ck(u),\bA(u)$ are defined as above but in terms of this modified $\Irr^*_CW$. However the two definitions are 
equivalent in view of the following result.

\proclaim{Proposition 0.2}Assume that $C$ is special. Let $E\in\Irr_CW$.

(a) We have $a_E\le e_C$.

(b) We have $a_E=e_C$ if and only if $E\in\cf$.
\endproclaim
This follows from \cite{\SUP, 10.9}. Note that (a) was stated without proof in \cite{\ORA, (13.1.2)} (the proof I
had in mind at the time of \cite{\ORA} was combinatorial).

\subhead 0.3\endsubhead
The following result is equivalent to a result stated without proof in \cite{\ORA, (13.1.3)}.

\proclaim{Theorem 0.4} Let $C$ be a special unipotent class of $G$, let $u\in C$ and let $\cf$ be the family that
contains $E_C$. Then we have canonically $\XX_\cf=M(\bA(u))$ so that the pairing $\{,\}$ on $\XX_\cf$ coincides 
with the pairing $\{,\}$ on $M(\bA(u))$. Hence $\cg_\cf$ can be taken to be $\bA(u)$.
\endproclaim
This is equivalent to the corresponding statement in the case where $G$ is adjoint, which reduces immediately to 
the case where $G$ is adjoint simple. It is then enough to prove the theorem for one $G$ in each isogeny class of
semisimple, almost simple algebraic groups; this will be done in \S3 after some combinatorial preliminaries in 
\S1,\S2. The proof uses the explicit description of the Springer correspondence: for type $A_n,G_2$ in 
\cite{\SP}; for type $B_n,C_n,D_n$ in \cite{\SH} (as an algorithm) and in \cite{\ICC} (by a closed formula); for 
type $F_4$ in \cite{\SHH}; for type $E_n$ in \cite{\AL},\cite{\SPA}.

An immediate consequence of (the proof of) Theorem 0.4 is the following result which answers a question of
R. Bezrukavnikov.

\proclaim{Corollary 0.5} In the setup of 0.4 let $E\in\Irr^*_CW$ and let $\cv_E$ be the corresponding 
$A(u)$-module viewed as an (irreducible) $\bA(u)$-module. The image of $E$ under the canonical imbedding
$\cf@>>>\XX_\cf=M(\bA(u))$ is represented by the pair $(1,\cv_E)\in M(\bA(u))$. Conversely, if $E\in\cf$ and the
image of $E$ under $\cf@>>>\XX_\cf=M(\bA(u))$ is represented by the pair $(1,\r)\in M(\bA(u))$ where $\r$ is an 
irreducible representation of $\bA(u)$, then $E\in\Irr^*_CW$ and $\r\cong\cv_E$.
\endproclaim

\subhead 0.6\endsubhead
Corollary 0.5 has the following interpretation. Let $Y$ be a (unipotent) character sheaf on $G$ whose restriction
to the regular semisimple elements is $\ne0$; assume that in the usual parametrization of 
unipotent character sheaves by $\sqc_{\cf'}\XX_{\cf'}$, $Y$ corresponds to $(1,\r)\in M(\bA(u))$ where $C$ is the
special unipotent class corresponding to a family $\cf$, $u\in C$ and $\r$ is an irreducible representation of 
$\bA(u)$. Then $Y|_C$ is (up to shift) the irreducible local system on $C$ defined by $\r$.

A parametrization of unipotent character sheaves on $G$ in terms of restrictions to various conjugacy classes of 
$G$ is outlined in \S4. 

\subhead 0.7. Notation\endsubhead
If $A,B$ are subsets of $\NN$ we denote by $A\dcup B$ the union of $A$ and $B$ regarded as a multiset (each 
element of $A\cap B$ appears twice). For any set $\cx$, we denote by $\cp(\cx)$ the set of subsets of $\cx$ viewed
as an $F_2$-vector space with sum given by the symmetric difference. If $\cx\ne\em$ we note that $\{\em,\cx\}$ is
a line in $\cp(\cx)$ and we set $\bcp(\cx)=\cp(\cx)/\{\em,\cx\}$, $\cp_{ev}(\cx)=\{L\in\cp(\cx);|L|=0\mod2\}$; let
$\bcp_{ev}(\cx)$ be the image of $\cp_{ev}(\cx)$ under the obvious map $\cp(\cx)@>>>\bcp(\cx)$ (thus 
$\bcp_{ev}(\cx)=\bcp(\cx)$ if $|\cx|$ is odd and $\bcp_{ev}(\cx)$ is a hyperplane in $\bcp(\cx)$ if $|\cx|$ is 
even). Now if $\cx\ne\em$, the assignment $L,L'\m|L\cap L'|\mod2$ defines a symplectic form on $\cp_{ev}(\cx)$
which induces a nondegenerate symplectic form $(,)$ on $\bcp_{ev}(\cx)$ via the obvious linear map 
$\cp_{ev}(\cx)@>>>\bcp_{ev}(\cx)$.

For $g\in G$ let $g_s$ (resp. $g_\o$) be the semisimple (resp. unipotent) part of $g$.

For $z\in(1/2)\ZZ$ we set $\lfloor z\rfloor=z$ if $z\in\ZZ$ and $\lfloor z\rfloor=z-(1/2)$ if $z\in\ZZ+(1/2)$.

{\it Erratum to \cite{\ORA}.} On page 86, line -6 delete: "$b'<b$" and on line -4 before "In the language..."
insert: "The array above is regarded as identical to the array obtained by interchanging its two rows."

On page 343, line -5, after "respect to $M$" insert: "and where the group $\cg_\cf$ defined in terms of $(u',M)$
is isomorphic to the group $\cg_\cf$ defined in terms of $(u,G)$".

{\it Erratum to \cite{\ICC}.} In the definition of $A_\a,B_\a$ in \cite{\ICC, 11.5}, the condition $I\in\a$ 
should be replaced by $I\in\a'$ and the condition $I\in\a'$ should be replaced by $I\in\a$.

\head 1. Combinatorics\endhead
\subhead 1.1\endsubhead
Let $N$ be an even integer $\ge0$. Let $a:=(a_0,a_1,a_2,\do,a_N)\in\NN^{N+1}$ be such that
$a_0\le a_1\le a_2\le\do\le a_N$, $a_0<a_2<a_4<\do$, $a_1<a_3<a_5<\do$. Let
$\cj=\{i\in[0,N];a_i\text{ appears exactly once in }a\}$. We have $\cj=\{i_0,i_1,\do,i_{2M}\}$ where $M\in\NN$ 
and $i_0<i_1<\do<i_{2M}$ satisfy $i_s=s\mod2$ for $s\in[0,2M]$. Hence for any $s\in[0,2M-1]$ we have 
$i_{s+1}=i_s+2m_s+1$ for some $m_s\in\NN$.
Let $\ce$ be the set of $b:=(b_0,b_1,b_2,\do,b_N)\in\NN^{N+1}$ such that $b_0<b_2<b_4<\do$, $b_1<b_3<b_5<\do$ and 
such that $[b]=[a]$ (we denote by $[b],[a]$ the multisets $\{b_0,b_1,\do,b_N\}$, $\{a_0,a_1,\do,a_N\}$). We have
$a\in\ce$. For $b\in\ce$ we set 
$$\hb=(\hb_0,\hb_1,\hb_2,\do,\hb_N)=(b_0,b_1+1,b_2+1,b_3+2,b_4+2,\do,b_{N-1}+(N/2),b_N+(N/2)).$$
Let $[\hb]$ be the multiset $\{\hb_0,\hb_1,\hb_2,\do,\hb_N\}$.

For $s\in\{1,3,\do,2M-1\}$ we define $a^{\{s\}}=(a^{\{s\}}_0,a^{\{s\}}_1,a^{\{s\}}_2,\do,a^{\{s\}}_N)\in\ce$ by 
$$\align&(a^{\{s\}}_{i_s},a^{\{s\}}_{i_s+1},a^{\{s\}}_{i_s+2},a^{\{s\}}_{i_s+3},\do,a^{\{s\}}_{i_s+2m_s},
a^{\{s\}}_{i_s+2m_s+1})\\&=(a_{i_s+1},a_{i_s},a_{i_s+3},a_{i_s+2},\do,a_{i_s+2m_s+1},a_{i_s+2m_s})\endalign$$
and $a^{\{s\}}_i=a_i$ if $i\in[0,N]-[i_s,i_{s+1}]$. More generally for $X\sub\{1,3,\do,2M-1\}$ we 
define $a^X=(a^X_0,a^X_1,a^X_2,\do,a^X_N)\in\ce$ by $a^X_i=a^{\{s\}}_i$ if $s\in X$, $i\in[i_s,i_{s+1}]$, and 
$a^X_i=a_i$ for all other $i\in[0,N]$. Note that $[\widehat{a^X}]=[\ha]$. Conversely, we have the following 
result.

\proclaim{Lemma 1.2} Let $b\in\ce$ be such that $[\hb]=[\ha]$. There exists $X\sub\{1,3,\do,2M-1\}$ such that
$b=a^X$.
\endproclaim
The proof is given in 1.3-1.5.

\subhead 1.3\endsubhead
We argue by induction on $M$. We have 
$$a=(y_1=y_1<y_2=y_2<\do<y_r=y_r<a_{i_0}<\do)$$ 
for some $r$. Since $[b]=[a]$, we must have
$$(b_0,b_2,b_4,\do)=(y_1,y_2,\do,y_r,\do),(b_1,b_3,b_5,\do)=(y_1,y_2,\do,y_r,\do).$$ 
Thus, 

(a) $b_i=a_i$ for $i<i_0$.
\nl
We have $a=(\do>a_{2M}>y'_1=y'_1<y'_2=y'_2<\do<y'_{r'}=y'_{r'})$ for some $r'$. Since $[b]=[a]$, we must have
$$(b_0,b_2,b_4,\do)=(\do,y'_1,y'_2,\do,y'_{r'}),(b_1,b_3,b_5,\do)=(\do,y'_1,y'_2,\do,y'_{r'}).$$ 
Thus, 

(b) $b_i=a_i$ for $i>i_{2M}$.
\nl
If $M=0$ we see that $b=a$ and there is nothing further to prove. In the rest of the proof we assume that 
$M\ge1$. 

\subhead 1.4\endsubhead
From 1.3 we see that
$$(a_0,a_1,a_2,\do,a_{i_{2M}})=(\do,a_{i_{2M-1}}<x_1=x_1<x_2=x_2<\do<x_q=x_q<a_{i_{2M}})$$
(for some $q$) has the same entries as $(b_0,b_1,b_2,\do,b_{i_{2M}}$ (in some order). Hence the pair
$$(\do,b_{i_{2M}-5},b_{i_{2M}-3},b_{i_{2M-1}}),(\do,b_{i_{2M}-4},b_{i_{2M}-2},b_{i_{2M}})$$
must have one of the following four forms.

$(\do,a_{i_{2M-1}},x_1,x_2,\do,x_q),(\do,x_1,x_2,\do,x_q,a_{i_{2M}})$,

$(\do,x_1,x_2,\do,x_q,a_{i_{2M}}),(\do,a_{i_{2M-1}},x_1,x_2,\do,x_q)$,

$(\do,x_1,x_2,\do,x_q),(\do,a_{i_{2M-1}},x_1,x_2,\do,x_q,a_{i_{2M}})$,

$(\do,a_{i_{2M-1}},x_1,x_2,\do,x_q,a_{i_{2M}}),(\do,x_1,x_2,\do,x_q)$.
\nl
Hence $(\do,b_{i_{2M}-2},b_{i_{2M}-1},b_{i_{2M}})$ must have one of the following four forms.

(I) $(\do,a_{i_{2M-1}},x_1,x_1,x_2,x_2,\do,x_q,x_q,a_{i_{2M}})$,

(II) $(\do,x_1,a_{i_{2M-1}},x_2,x_1,x_3,x_2,\do,x_q,x_{q-1},a_{i_{2M}},x_q)$,

(III) $(\do,a_{i_{2M-1}},z,x_1,x_1,x_2,x_2,\do,x_q,x_q,a_{i_{2M}})$,

(IV) $(\do,a_{i_{2M-1}},z',x_1,z'',x_2,x_1,x_3,x_2,\do,x_q,x_{q-1},a_{i_{2M}},x_q)$,
\nl
where $a_{i_{2M-1}}>z$,  $a_{i_{2M-1}}>z''\ge z'$ and all entries in $\do$ are $<a_{i_{2M-1}}$. Correspondingly,
$(\do,\hb_{i_{2M}-2},\hb_{i_{2M}-1},\hb_{i_{2M}})$ must have one of the following four forms.

(I) $(\do,a_{i_{2M-1}}+h-q,x_1+h-q,x_1+h-q+1,x_2+h-q+1,x_2+h-q+2,\do,x_q+h-1,x_q+h,a_{i_{2M}}+h)$,

(II) $(\do,x_1+h-q,a_{i_{2M-1}}+h-q,x_2+h-q+1,x_1+h-q+1,x_3+h-q+2,x_2+h-q+1,\do,x_q+h-1,x_{q-1}+h-1,
a_{i_{2M}}+h,x_q+h)$,

(III) $(\do,a_{i_{2M-1}}+h-q-1,z+h-q,x_1+h-q,x_1+h-q+1,x_2+h-q+1,x_2+h-q+2,\do,x_q+h-1,x_q+h,a_{i_{2M}}+h)$,

(IV) $(\do,a_{i_{2M-1}}+h-q-1,z'+h-q-1,x_1+h-q,z''+h-q,x_2+h-q+1,x_1+h-q+1,x_3+h-q+2,x_2+h-q+1,\do,x_q+h-1,
x_{q-1}+h-1,a_{i_{2M}}+h,x_q+h)$
\nl
where $h=i_{2M}/2$ and in case (III) and (IV), $a_{i_{2M-1}}+h-q$ is not an entry of 
$(\do,\hb_{i_{2M}-2},\hb_{i_{2M}-1},\hb_{i_{2M}})$.

Since $(\do,\ha_{i_{2M}-2},\ha_{i_{2M}-1},\ha_{i_{2M}})$ is given by (I) we see that 
$a_{i_{2M-1}}+h-q$ is an entry of $(\do,\ha_{i_{2M}-2},\ha_{i_{2M}-1},\ha_{i_{2M}})$. Using 1.3(b) we see that 

$\{\do,\ha_{i_{2M}-2},\ha_{i_{2M}-1},\ha_{i_{2M}}\}=(\do,b_{i_{2M}-2},b_{i_{2M}-1},b_{i_{2M}})$
\nl
as multisets. We see that cases (III) and (IV) cannot arise. Hence we must be in case (I) or (II). Thus 

{\it we have either
$$\align&(b_{i_{2M-1}},b_{i_{2M-1}+1},\do,b_{i_{2M}-2},b_{i_{2M}-1},b_{i_{2M}})\\&=
(a_{i_{2M-1}},a_{i_{2M-1}+1},\do,a_{i_{2M}-2},a_{i_{2M}-1},a_{i_{2M}})\tag a\endalign$$
or
$$\align&(b_{i_{2M-1}},b_{i_{2M-1}+1},\do,b_{i_{2M}-2},b_{i_{2M}-1},b_{i_{2M}})\\&=
(a_{i_{2M-1}+1},a_{i_{2M-1}},a_{i_{2M-1}+3},a_{i_{2M-1}+2},\do,a_{i_{2M}},a_{i_{2M}-1}).\tag b\endalign$$ }

\subhead 1.5\endsubhead
Let $a'=(a_0,a_1,a_2,\do,a_{i_{2M-1}-1})$, $b'=(b_0,b_1,b_2,\do,b_{i_{2M-1}-1})$,

$\ha'=(a_0,a_1+1,a_2+1,a_3+2,a_4+2,\do,a_{i_{2M-1}-1}+(i_{2M-1}-1)/2)$,

$\hb'=(b_0,b_1+1,b_2+1,b_3+2,b_4+2,\do,b_{i_{2M-1}-1}+(i_{2M-1}-1)/2)$,
\nl
From $[\hb]=[\ha]$ and 1.3(b), 1.4(a),(b) we see that the multiset formed by the entries of $\ha'$ coincides with
the multiset formed by the entries of $\hb'$. Using the induction hypothesis we see that there exists 
$X'\sub\{1,3,\do,2M-3\}$ such that $b'=a'{}^{X'}$ where $a'{}^{X'}$ is defined in terms of $a',X'$ in the same 
way as $a^X$ was defined (see 1.1) in terms of $a,X$. We set $X=X'$ if we are in case 1.4(a) and 
$X=X'\cup\{2M-1\}$ if we are in case 1.4(b). Then we have $b=a^X$ (see 1.4(a),(b)), as required. This completes 
the proof of Lemma 1.2.

\subhead 1.6\endsubhead
We shall use the notation of 1.1. Let $\fT$ be the set of all unordered pairs $(\fA,\fB)$ where $\fA,\fB$ are
subsets of $\{0,1,2,\do\}$ and $\fA\dcup\fB=(a_0,a_1,a_2,\do,a_N)$ as multisets. For example, setting 
$\fA_\em=(a_0,a_2,a_4,\do,a_N)$, $\fB_\em=(a_1,a_3,\do,a_{N-1})$, we have $(\fA_\em,\fB_\em)\in\fT$. For any 
subset $\fa$ of $\cj$ we consider
$$\fA_{\fa}=((\cj-\fa)\cap\fA_\em)\cup(\fa\cap\fB_\em)\cup(\fA_\em\cap\fB_\em),$$
$$\fB_{\fa}=((\cj-\fa)\cap\fB_\em)\cup(\fa\cap\fA_\em)\cup(\fA_\em\cap\fB_\em).$$
Then $(\fA_{\fa},\fB_{\fa})\in\fT$ and the map $\fa\m(\fA_{\fa},\fB_{\fa})$ induces a bijection 
$\bcp(\cj)\lra\fT$. (Note that if $\fa=\em$ then $(\fA_{\fa},\fB_{\fa})$ agrees with the earlier definition of 
$(\fA_\em,\fB_\em)$.) 

Let $\fT'$ be the set of all $(\fA,\fB)\in\fT$ such that $|\fA|=|\fA_\em|,|\fB|=|\fB_\em|$.

Let $\cp(\cj)_0$ be the subspace of $\cp_{ev}(\cj)$ spanned by the $2$-element subsets
$$\{a_{i_0},a_{i_1}\},\{a_{i_2},a_{i_3}\},\do,\{a_{i_{2M-2}},a_{i_{2M-1}}\}$$ 
of $\cj$. Let $\cp(\cj)_1$ be the subspace of $\cp_{ev}(\cj)$ spanned by the $2$-element subsets
$$\{a_{i_1},a_{i_2}\},\{a_{i_3},a_{i_4}\},\do,\{a_{i_{2M-1}},a_{i_{2M}}\}$$ 
of $\cj$. 

Let $\bcp(\cj)_0$ (resp. $\bcp(\cj)_1$) be the image of $\cp(\cj)_0$ (resp. $\cp(\cj)_1$) under the obvious map 
$\cp(\cj)@>>>\bcp(\cj)$. Note that 

(a) {\it $\bcp(\cj)_0$ and $\bcp(\cj)_1$ are opposed Lagrangian subspaces of the symplectic vector space 
$\bcp(\cj),(,)$, (see 0.7); hence $(,)$ defines an identification $\bcp(\cj)_0=\bcp(\cj)_1^*$}
\nl
where $\bcp(\cj)_1^*$ is the vector space dual to $\bcp(\cj)_1$.

Let $\fT_0$ (resp. $\fT_1$) be the subset of $\fT$ corresponding to $\bcp(\cj)_0$ (resp. $\bcp(\cj)_1$)
under the bijection $\bcp(\cj)\lra\fT$. Note that  $\fT_0\sub\fT'$, $\fT_1\sub\fT'$, $|\fT_0|=|\fT_1|=2^M$.

For any $X\sub\{1,3,\do,2M-1\}$ we set $\fa_X=\cup_{s\in X}\{a_{i_s},a_{i_{s+1}}\}\in\cp(\cj)$. Then 
$(\fA_{\fa_X},\fB_{\fa_X})\in\fT_1$ is related to $a^X$ in 1.1 as follows:

$\fA_{\fa_X}=\{a^X_0,a^X_2,a^X_4,\do,a^X_N\}$, $\fB_{\fa_X}=\{a^X_1,a^X_3,\do,a^X_{N-1}\}$.

\subhead 1.7\endsubhead
We shall use the notation of 1.1. Let $T$ be the set of all ordered pairs $(A,B)$ where $A$ is a subset of 
$\{0,1,2,\do\}$, $B$ is a subset of $\{1,2,3,\do\}$, $A$ contains no consecutive integers, $B$ contains no 
consecutive integers, and $A\dcup B=(\ha_0,\ha_1,\ha_2,\do,\ha_N)$ as multisets. For example, setting 
$A_\em=(\ha_0,\ha_2,\ha_4,\do,\ha_N)$, $B_\em=(\ha_1,\ha_3,\do,\ha_{N-1})$, we have $(A_\em,B_\em)\in T$.

For any $(A,B)\in T$ we define $(A^-,B^-)$ as follows: $A^-$ consists of $x_0<x_1-1<x_2-2<\do<x_p-p$ where 
$x_0<x_1<\do<x_p$ are the elements of $A$; $B^-$ consists of $y_1-1<y_2-2<\do<y_q-q$ where $y_1<y_2<\do<y_q$ are 
the elements of $B$. 

We can enumerate the elements of $T$ as in \cite{\ICC, 11.5}. Let $J$ be the set of all $c\in\NN$ such that $c$ 
appears exactly once in the sequence 
$$(\ha_0,\ha_1,\ha_2,\do,\ha_N)=(a_0,a_1+1,a_2+1,a_3+2,a_4+2,\do,a_{N-1}+(N/2),a_N+(N/2)).$$
A nonempty subset $I$ of $J$ is said to be an interval if it is of the form $\{i,i+1,i+2,\do,j\}$ with 
$i-1\n J,j+1\n J$ and with $i\ne0$. Let $\ci$ be the set of intervals of $J$. For any $s\in\{1,3,\do,2M-1\}$, the
set $I_s:=\{\ha_{i_s},\ha_{i_s+1},\ha_{i_s+2},\do,\ha_{i_s+2m_s+1}\}$ is either a single interval or a union of 
intervals $I_s^1\sqc I_s^2\sqc\do\sqc I_s^{t_s}$ ($t_s\ge2$) where $\ha_{i_s}\in I_s^1$, 
$\ha_{i_s+2m_s+1}\in I_s^{t_s}$, $|I_s^1|,|I_s^{t_s}|$ are odd, $|I_s^h|$ are even for $h\in[2,t_s-1]$ and any 
element in $I_s^e$ is $<$ than any element in $I_s^{e'}$ for $e<e'$. Let $\ci_s$ be the set of all $I\in\ci$ such
that $I\sub I_s$. We have a partition $\ci=\sqc_{s\in\{1,3,\do,2M-1\}}\ci_s$. Let $H$ be the set of elements of 
$c\in J$ such that $c<a_{i_1}$ (that is such that $c$ does not belong to any interval). For any subset 
$\a\sub\ci$ we consider 
$$A_\a=\cup_{I\in\ci-\a}(I\cap A_\em)\cup\cup_{I\in\a}(I\cap B_\em)\cup(H\cap A_\em)\cup(A_\em\cap B_\em),$$
$$B_\a=\cup_{I\in\ci-\a}(I\cap B_\em)\cup\cup_{I\in\a}(I\cap A_\em)\cup(H\cap B_\em)\cup(A_\em\cap B_\em).$$
Then $(A_\a,B_\a)\in T$ and the map $\a\m(A_\a,B_\a)$ is a bijection $\cp(\ci)\lra T$. (Note that if $\a=\em$ 
then $(A_\a,B_\a)$ agrees with the earlier definition of $(A_\em,B_\em)$.) 

Let $T'=\{(A,B)\in T;|A|=|A_\em|,|B|=|B_\em|\}$, $T_1=\{(A,B)\in T';A^-\dcup B^-=A_\em^-\dcup B_\em^-\}$. Let 
$\cp(\ci)'$ (resp. $\cp(\ci)_1$) be the subset of $\cp(\ci)$ corresponding to $T'$ (resp. $T_1$) under the 
bijection $\cp(\ci)\lra T$.

Now let $X$ be a subset of $\{1,3,\do,2M-1\}$. Let $\a_X=\cup_{s\in X}\ci_s\in\cp(\ci)$. From the definitions we 
see that 

(a) $A_{\a_X}^-=\fA_{\fa_X}$, $B_{\a_X}^-=\fB_{\fa_X}$
\nl
(notation of 1.6). In particular we have $(A_{\a_X},B_{\a_X})\in T_1$. Thus $|T_1|\ge2^M$. Using Lemma 1.2 we see
that

(b) {\it $|T_1|=2^M$ and $T_1$ consists of the pairs $(A_{\a_X},B_{\a_X})$ with $X\sub\{1,3,\do,2M-1\}$.}
\nl
Using (a),(b) we deduce:

(c) {\it The map $T_1@>>>\fT_1$ given by $(A,B)\m(A^-,B^-)$ is a bijection.}

\head 2. Combinatorics (continued)\endhead
\subhead 2.1\endsubhead
Let $N\in\NN$. Let 
$$a:=(a_0,a_1,a_2,\do,a_N)\in\NN^{N+1}$$ be such that 
$a_0\le a_1\le a_2\le\do\le a_N$, $a_0<a_2<a_4<\do$, $a_1<a_3<a_5<\do$ and such that the set 
$\cj:=\{i\in[0,N];a_i\text{ appears exactly once in }a\}$ is nonempty. Now $\cj$ consists of $\mu+1$ elements
$i_0<i_1<\do<i_\mu$ where $\mu\in\NN$, $\mu=N\mod2$. We have $i_s=s\mod2$ for $s\in[0,\mu]$. Hence for any 
$s\in[0,\mu-1]$ we have $i_{s+1}=i_s+2m_s+1$ for some $m_s\in\NN$. Let $\ce$ be the set of 
$b:=(b_0,b_1,b_2,\do,b_N)\in\NN^{N+1}$ such that $b_0<b_2<b_4<\do$, $b_1<b_3<b_5<\do$ and such that $[b]=[a]$ (we
denote by $[b],[a]$ the multisets $\{b_0,b_1,\do,b_N\}$, $\{a_0,a_1,\do,a_N\}$). We have $a\in\ce$. For $b\in\ce$
we set 
$$\ocb=(\ocb_0,\ocb_1,\ocb_2,\do,\ocb_N)=(b_0,b_1,b_2+1,b_3+1,b_4+2,b_5+2,\do)\in\NN^{N+1}.$$
Let $[\ocb]$ be the multiset $\{\ocb_0,\ocb_1,\ocb_2,\do,\ocb_N\}$. For any $s\in[0,\mu-1]\in2\NN$ we define 
$a^{\{s\}}=(a^{\{s\}}_0,a^{\{s\}}_1,a^{\{s\}}_2,\do,a^{\{s\}}_N)\in\ce$ by 
$$\align&(a^{\{s\}}_{i_s},a^{\{s\}}_{i_s+1},a^{\{s\}}_{i_s+2},a^{\{s\}}_{i_s+3},\do,a^{\{s\}}_{i_s+2m_s},
a^{\{s\}}_{i_s+2m_s+1})\\&=(a_{i_s+1},a_{i_s},a_{i_s+3},a_{i_s+2},\do,a_{i_s+2m_s+1},a_{i_s+2m_s})\endalign$$
and $a^{\{s\}}_i=a_i$ if $i\in[0,N]-[i_s,i_{s+1}]$. More generally for a subset $X$ of $[0,\mu-1]\cap2\NN$ we 
define $a^X=(a^X_0,a^X_1,a^X_2,\do,a^X_N)\in\ce$ by $a^X_i=a^{\{s\}}_i$ if $s\in X$, $i\in[i_s,i_{s+1}]$, and 
$a^X_i=a_i$ for all other $i\in[0,N]$. Note that $[\oca^X]=[\oca]$. Conversely, we have the following result.

\proclaim{Lemma 2.2} Let $b\in\ce$ be such that $[\ocb]=[\oca]$. Then there exists $X\sub[0,\mu-1]\cap2\NN$ such 
that $b=a^X$.
\endproclaim
The proof is given in 2.3-2.5.

\subhead 2.3\endsubhead
We argue by induction on $\mu$. By the argument in 1.3 we have

(a) $b_i=a_i$ for $i<i_0$,

(b) $b_i=a_i$ for $i>i_\mu$.
\nl
If $\mu=0$ we see that $b=a$ and there is nothing further to prove. In the rest of the proof we assume that 
$\mu\ge1$. 

\subhead 2.4\endsubhead
From 2.3 we see that $(a_{i_0},a_{i_0+1},\do,a_N)=(a_{i_0}<x_1=x_1<x_2=x_2<\do<x_p=x_p<a_{i_1}<\do)$ (for some 
$p$) has the same entries as $(b_{i_0},b_{i_0+1},\do,b_N)$ (in some order). Hence the pair
$(b_{i_0},b_{i_0+2},b_{i_0+4},\do),(b_{i_0+1},b_{i_0+3},b_{i_0+5},\do)$ must have one of the following four forms.

$(a_{i_0},x_1,x_2,\do,x_p,\do),(x_1,x_2,\do,x_p,a_{i_1},\do)$,

$(x_1,x_2,\do,x_p,a_{i_1},\do),(a_{i_0},x_1,x_2,\do,x_p,\do)$,

$(a_{i_0},x_1,x_2,\do,x_p,a_{i_1},\do),(x_1,x_2,\do,x_p,\do)$,

$(x_1,x_2,\do,x_p,\do),(a_{i_0},x_1,x_2,\do,x_p,a_{i_1},\do)$.
\nl
Hence $(b_{i_0},b_{i_0+1},b_{i_0+2},\do,b_N)$ must have one of the following four forms.

(I) $(a_{i_0},x_1,x_1,x_2,x_2,\do,x_p,x_p,a_{i_1},\do)$,

(II) $(x_1,a_{i_0},x_2,x_1,x_3,x_2,\do,x_p,x_{p-1},a_{i_1},x_p,\do)$,

(III) $(a_{i_0},x_1,x_1,x_2,x_2,\do,x_p,x_p,z,a_{i_1},\do)$,

(IV) $(x_1,a_{i_0},x_2,x_1,x_3,x_2,\do,x_p,x_{p-1},z',x_p,z'',a_{i_1},\do)$
\nl
where $a_{i_1}<z$, $a_{i_1}<z'\le z''$ and all entries in $\do$ are $>a_{i_1}$.
Correspondingly, $(\ocb_{i_0},\ocb_{i_0+1},\ocb_{i_0+2},\do,\ocb_N)$ must have one of the following four forms.

(I) $(a_{i_0}+h,x_1+h,x_1+h+1, x_2+h+1,x_2+h+2,\do,x_p+h+p-1,x_p+h+p,a_{i_1}+h+p,\do)$,

(II) $(x_1+h,a_{i_0}+h,x_2+h+1,x_1+h+1,x_3+h+2,x_2+h+2,\do,x_p+h+p-1,x_{p-1}+h+p-1,a_{i_1}+h+p,x_p+h+p,\do)$,

(III) $(a_{i_0}+h,x_1+h,x_1+h+1,x_2+h+1,x_2+h+2,\do,x_p+h+p-1,x_p+h+p,z+p,a_{i_1}+h+p+1,\do)$,

$$\align&(x_1+h,a_{i_0}+h,x_2+h+1,x_1+h+1,x_3+h+2,x_2+h+2,\do,x_p+h+p-1,\\&
x_{p-1}+h+p-1,z'+h+p,x_p+h+p,z''+h+p+1,a_{i_1}+h+p+1,\do) \tag IV\endalign$$

where $h=i_0/2$ and in case (III) and (IV) $a_{i_1}+h+p$ is not an entry of 
$(\ocb_{i_0},\ocb_{i_0+1},\ocb_{i_0+2},\do)$.

Since $(\oca_{i_0},\oca_{i_0+1},\oca_{i_0+2},\do)$ is given by (I) we see that $a_{i_1}+h+p$ is an entry of 
$(\oca_{i_0},\oca_{i_0+1},\oca_{i_0+2},\do)$. Using 2.3 we see that 
$$\{\oca_{i_0},\oca_{i_0+1},\oca_{i_0+2},\do\}=\{\ocb_{i_0},\ocb_{i_0+1},\ocb_{i_0+2},\do\}$$
as multisets. We see that cases (III) and (IV) cannot arise. Hence we must be in case (I) or (II). Thus 

{\it we have either
$$(b_{i_0},b_{i_0+1},b_{i_0+2},\do,b_{i_1})=(a_{i_0},a_{i_0+1},a_{i_0+2},\do,a_{i_1})\tag a$$
or 
$$(b_{i_0},b_{i_0+1},b_{i_0+2},\do,b_{i_1})=(a_{i_0+1},a_{i_0},a_{i_0+3},a_{i_0+2},\do,a_{i_1},a_{i_1-1}).\tag b
$$.}
From 2.3 and (a),(b) we see that if $\mu=1$ then Lemma 2.2 holds. Thus in the rest of the proof we can assume 
that $\mu\ge2$.

\subhead 2.5\endsubhead
Let $a'=(a_{i_1+1},a_{i_1+2},\do,a_N)$, $b'=(b_{i_1+1},b_{i_1+2},\do,b_N)$,
$$\oca'=(a_{i_1+1},a_{i_1+2},a_{i_1+3}+1,a_{i_1+4}+1, a_{i_1+5}+2,a_{i_1+6}+2,\do),$$ 
$$\ocb'=(b_{i_1+1},b_{i_1+2},b_{i_1+3}+1,b_{i_1+4}+1, b_{i_1+5}+2,b_{i_1+6}+2,\do).$$
From $[\ocb]=[\oca]$ and 2.3(a),2.4(a),(b) we see that the multiset formed by the entries of $\oca'$ coincides 
with the multiset formed by the entries of $\ocb'$. Using the induction hypothesis we see that there exists 
$X'\sub[2,\mu-1]\cap2\NN$ such that $b'=a'{}^{X'}$ where $a'{}^{X'}$ is defined in terms of $a',X'$ in the same 
way as $a^X$ (see 2.1) was defined in terms of $a,X$. We set $X=X'$ if we are in case 2.4(a) and $X=\{0\}\cup X'$
if we are in case 2.4(b). Then we have $b=a^X$ (see 2.4(a),(b)), as required. This completes the proof of Lemma 
2.2.

\subhead 2.6\endsubhead
We shall use the notation of 2.1. Let $\fT$ be the set of all unordered pairs $(\fA,\fB)$ where $\fA,\fB$ are
subsets of $\{0,1,2,\do\}$ and $\fA\dcup\fB=(a_0,a_1,a_2,\do,a_N)$ as multisets. For example, setting 
$\fA_\em=\{a_i;i\in[0,N]\cap2\NN\}$, $\fB_\em=\{a_i;i\in[0,N]\cap(2\NN+1)\}$, we have $(\fA_\em,\fB_\em)\in\fT$.
For any subset $\fa$ of $\cj$ we consider
$$\fA_{\fa}=((\cj-\fa)\cap\fA_\em)\cup(\fa\cap\fB_\em)\cup(\fA_\em\cap\fB_\em),$$
$$\fB_{\fa}=((\cj-\fa)\cap\fB_\em)\cup(\fa\cap\fA_\em)\cup(\fA_\em\cap\fB_\em).$$
Then $(\fA_{\fa},\fB_{\fa})=(\fA_{\cj-\fa},\fA_{\cj-\fa})\in\fT$ and the map $\fa\m(\fA_{\fa},\fB_{\fa})$ induces
a bijection $\bcp(\cj)\lra\fT$. (Note that if $\fa=\em$ then $(\fA_{\fa},\fB_{\fa})$ agrees with the earlier 
definition of $(\fA_\em,\fB_\em)$.) 

Let $\fT'$ be the set of all $(\fA,\fB)\in\fT$ such that $|\fA|=|\fA_\em|,|\fB|=|\fB_\em|$. 
Let $\cp(\cj)_1$ be the subspace of $\cp(\cj)$ spanned by the following $2$-element subsets of $\cj$:
$$\{a_{i_1},a_{i_2}\},\{a_{i_3},a_{i_4}\},\do,\{a_{i_{\mu-2}},a_{i_{\mu-1}}\}\text{(if $N$ is odd)}$$
or
$$\{a_{i_1},a_{i_2}\},\{a_{i_3},a_{i_4}\},\do,\{a_{i_{\mu-1}},a_{i_{\mu}}\}\text{(if $N$ is even)}.$$
Let $\cp(\cj)_0$ be 
the subspace of $\cp(\cj)$ spanned by the following $2$-element subsets of $\cj$:
$$\{a_{i_0},a_{i_1}\},\{a_{i_2},a_{i_3}\},\do,\{a_{i_{\mu-1}},a_{i_\mu}\}\text{(if $N$ is odd)}$$
or
$$\{a_{i_0},a_{i_1}\},\{a_{i_2},a_{i_3}\},\do,\{a_{i_{\mu-2}},a_{i_{\mu-1}}\}\text{(if $N$ is even)}.$$
Let $\bcp(\cj)_0$ (resp. $\bcp(\cj)_1$) be the image of $\cp(\cj)_0$ (resp. $\cp(\cj)_1$) under the obvious map 
$\cp(\cj)@>>>\bcp(\cj)$. 

Note that 

(a) {\it $\bcp(\cj)_0$ and $\bcp(\cj)_1$ are opposed Lagrangian subspaces of the symplectic vector space 
$\bcp_{ev}(\cj),(,)$, (see 0.7); hence $(,)$ defines an identification $\bcp(\cj)_1=\bcp(\cj)_0^*$}
\nl
where $\bcp(\cj)_0^*$ is the vector space dual to $\bcp(\cj)_0$.

Let $\fT_0$ (resp. $\fT_1$) be the subset of $\fT$ corresponding to $\bcp(\cj)_0$ (resp. $\bcp(\cj)_1$) under the
bijection $\bcp(\cj)\lra\fT$. Note that $\fT_0\sub\fT',\fT_1\sub\fT'$, $|\fT_0|=|\fT_1|=2^{\lfloor\mu/2\rfloor}$.

For any $X\sub[0,\mu-1]\cap2\NN$ we set $\fa_X=\cup_{s\in X}\{a_{i_s},a_{i_{s+1}}\}\in\cp(\cj)$. Then 
$(\fA_{\fa_X},\fB_{\fa_X})$ is related to $a^X$ in 2.1 as follows:
$$\fA_{\fa_X}=\{a^X_i;i\in[0,N]\cap2\NN\}, \fB_{\fa_X}=\{a^X_i;i\in[0,N]\cap(2\NN+1)\}.$$
     
\subhead 2.7\endsubhead
We shall use the notation of 2.1. Let $T$ be the set of all unordered pairs $(A,B)$ where $A$ is a subset of 
$\{0,1,2,\do\}$, $B$ is a subset of $\{1,2,3,\do\}$, $A$ contains no consecutive integers, $B$ contains no 
consecutive integers, and $A\dcup B=(\oca_0,\oca_1,\oca_2,\do,\oca_N)$ as multisets. For example, setting
$$A_\em=\{\oca_i;i\in[0,N]\cap2\NN\}, B_\em=(\oca_i;i\in[0,N]\cap(2\NN+1)\},$$
we have $(A_\em,B_\em)\in T$.

For any $(A,B)\in T$ we define $(A^-,B^-)$ as follows: $A^-$ consists of $x_1<x_2-1<x_3-2<\do<x_p-p+1$ where 
$x_1<x_2<\do<x_p$ are the elements of $A$; $B^-$ consists of $y_1<y_2-1<\do<y_q-q+1$ where $y_1<y_2<\do<y_q$ are 
the elements of $B$. 

We can enumerate the elements of $T$ as in \cite{\ICC, 11.5}. Let $J$ be the set of all $c\in\NN$ such that $c$ 
appears exactly once in the sequence 
$$(\oca_0,\oca_1,\oca_2,\do,\oca_N)=(a_0,a_1,a_2+1,a_3+1,a_4+2,a_5+2,\do).$$
A nonempty subset $I$ of $J$ is said to be an interval if it is of the form $\{i,i+1,i+2,\do,j\}$ with 
$i-1\n J,j+1\n J$. Let $\ci$ be the set of intervals of $J$. For any $s\in[0,\mu-1]\cap2\NN$, the
set $I_s:=\{\oca_{i_s},\oca_{i_s+1},\oca_{i_s+2},\do,\oca_{i_s+2m_s+1}\}$ is either a single interval or a union 
of intervals $I_s^1\sqc I_s^2\sqc\do\sqc I_s^{t_s}$ ($t_s\ge2$) where $\oca_{i_s}\in I_s^1$, 
$\oca_{i_s+2m_s+1}\in I_s^{t_s}$, $|I_s^1|,|I_s^{t_s}|$ are odd, $|I_s^h|$ are even for $h\in[2,t_s-1]$ and any 
element in $I_s^e$ is $<$ than any element in $I_s^{e'}$ for $e<e'$. Let $\ci_s$ be the set of all $I\in\ci$ such
that $I\sub I_s$. We have a partition $\ci=\sqc_{s\in[0,\mu-1]\cap2\NN}\ci_s$. For any subset $\a\sub\ci$ we 
consider 
$$A_\a=\cup_{I\in\ci-\a}(I\cap A_\em)\cup\cup_{I\in\a}(I\cap B_\em)\cup(A_\em\cap B_\em),$$
$$B_\a=\cup_{I\in\ci-\a}(I\cap B_\em)\cup\cup_{I\in\a}(I\cap A_\em)\cup(A_\em\cap B_\em).$$
Then $(A_\a,B_\a)\in T$ and the map $\a\m(A_\a,B_\a)$ is a bijection $\bcp(\ci)\lra T$.
(Note that if $\a=\em$ then $(A_\a,B_\a)$ agrees with the earlier definition of $(A_\em,B_\em)$.) 

Let $T'=\{(A,B)\in T;|A|=|A_\em|,|B|=|B_\em|\}$, $T_1=\{(A,B)\in T';A^-\dcup B^-=A_\em^-\dcup B_\em^-\}$. Let 
$\bcp(\ci)'$ (resp. $\bcp(\ci)_1$) be the subset of $\bcp(\ci)$ corresponding to $T'$ (resp. $T_1$) under the 
bijection $\bcp(\ci)\lra T$.

Now let $X$ be a subset of $[0,\mu-1]\cap2\NN$. Let $\a_X=\cup_{s\in X}\ci_s\in\cp(\ci)$. From the definitions we 
see that 

(a) $A_{\a_X}^-=\fA_{\fa_X}$, $B_{\a_X}^-=\fB_{\fa_X}$
\nl
(notation of 2.6). In particular we have $(A_{\a_X},B_{\a_X})\in T_1$. Thus $|T_1|\ge2^{\lfloor\mu/2\rfloor}$.
Using Lemma 2.2 we see that

(b) {\it $|T_1|=2^{\lfloor\mu/2\rfloor}$ and $T_1$ consists of the pairs $(A_{\a_X},B_{\a_X})$ with 
$X\sub[0,\mu-1]\cap2\NN$.}
\nl
Using (a),(b) we deduce:

(c) {\it The map $T_1@>>>\fT_1$ given by $(A,B)\m(A^-,B^-)$ is a bijection.}

\head 3. Proof of Theorem 0.4 and of Corollary 0.5\endhead
\subhead 3.1\endsubhead
If $G$ is simple adjoint of type $A_n$, $n\ge1$, then 0.4 and 0.5 are obvious: we have $A(u)=\{1\}$, 
$\bA(u)=\{1\}$.

\subhead 3.2\endsubhead
Assume that $G=Sp_{2n}(\kk)$ where $n\ge2$. Let $N$ be a sufficiently large even integer. Now 
$u:\kk^{2n}@>>>\kk^{2n}$ has $i_e$ Jordan blocks of size $e$ ($e=1,2,3,\do$). Here $i_1,i_3,i_5,\do$ are even. 
Let $\D=\{e\in\{2,4,6,\do\};i_e\ge1\}$. Then $A(u)$ can be identified in the standard way with $\cp(\D)$. Hence 
the group of characters $\hA(u)$ of $A(u)$ (which may be canonically identified with the $F_2$-vector space dual 
to $\cp(\D)$) may be also canonically identified with $\cp(\D)$ itself (so that the basis given by the one element
subsets of $\D$ is self-dual).

To the partition $1i_1+2i_2+3i_3+\do$ of $2n$ we associate a pair $(A,B)$ as in \cite{\ICC, 11.6} (with $N,2m$ 
replaced by $2n,N$). We have $A=(\ha_0,\ha_2,\ha_4,\do,\ha_N)$, $B=(\ha_1,\ha_3,\do,\ha_{N-1})$, where 
$\ha_0\le\ha_1\le\ha_2\le\do\le\ha_N$ is obtained from a sequence $a_0\le a_1\le a_2\le\do\le a_N$ as in 1.1. 
(Here we use that $C$ is special.) Now the definitions and results in \S1 are applicable. As in \cite{\ORA, 4.5} 
the family $\cf$ is in canonical bijection with $\fT'$ in 1.6.

We arrange the intervals in $\ci$ in increasing order $I_{(1)},I_{(2)},\do,I_{(f)}$ (any element in $I_{(1)}$ is 
smaller than any element in $I_{(2)}$, etc.). We arrange the elements of $\D$ in increasing order 
$e_1<e_2<\do<e_{f'}$; then $f=f'$ and we have a bijection $\ci\lra\D$, $I_{(h)}\lra e_h$; moreover we have 
$|I_{(h)}|=i_{e_h}$ for $h\in[1,f]$; see \cite{\ICC, 11.6}. Using this bijection we see that $A(u)$ and $\hA(u)$ 
are identified with the $F_2$-vector space $\cp(\ci)$ with basis given by the one element subsets of $\ci$. Let 
$\p:\cp(\ci)@>>>\cp(\ci)_1^*$ (with $\cp(\ci)_1^*$ as in 1.7(c)) be the (surjective) $F_2$-linear map which to 
$X\sub\ci$ associates the linear form $L\m|X\cap L|\mod2$ on $\cp(\ci)_1$. We will show that 

(a) $\ker\p=\ck(u)$ ($\ck(u)$ as in 0.1).
\nl
We identify $\Irr_CW$ with $T'$ (see 1.7) via the restriction of the bijection in \cite{\ICC, (12.2.4)} (we also
use the description of the Springer correspondence in \cite{\ICC, 12.3}). Under this identification the subset 
$\Irr^*_CW$ of $\Irr_CW$ becomes the subset $T_1$ (see 1.7) of $T'$. Via the identification $\cp(\ci)'\lra T'$ in
1.7 and $\hA(u)\lra\cp(\ci)$ (see above), the map $E\m\cv_E$ from $T'$ to $\hA(u)$ becomes the obvious imbedding 
$\cp(\ci)'@>>>\cp(\ci)$ (we use again \cite{\ICC, 12.3}). By definition, $\ck(u)$ is the set of all 
$X\in\cp(\ci)$ such that for any $L\in\cp(\ci)_1$ we have $|X\cap L|=0\mod2$. Thus, (a) holds.

Using (a) we have canonically $\bA(u)=\cp(\ci)_1^*$ via $\p$. We define an $F_2$-linear map 
$\cp(\ci)_1@>>>\bcp(\cj)_1$ (see 1.6) by $I_s\m\{a_{i_s},a_{i_{s+1}}\}$ for $s\in\{1,3,\do,2M-1\}$ ($I_s$ as in 
1.7). This is an isomorphism; it corresponds to the bijection 1.7(c) under the identification
$T_1\lra\cp(\ci)_1$ in 1.7 and the identification $\fT_1\lra\bcp(\cj)_1$ in 1.6. Hence we can identify 
$\cp(\ci)_1^*$ with $\bcp(\cj)_1^*$ and with $\bcp(\cj)_0$ (see 1.6(a)). We obtain an identification 
$\bA(u)=\bcp(\cj)_0$.

By \cite{\ORA, 4.5} we have $\XX_\cf=\bcp(\cj)$. Using 1.6(a) we see that \lb 
$\bcp(\cj)=M(\bcp(\cj)_0)=M(\bA(u))$ canonically so that 0.4 holds in our case. From the arguments above we see 
that in our case 0.5 follows from 1.7(c).

\subhead 3.3\endsubhead
Assume that $G=SO_n(\kk)$ where $n\ge7$. Let $N$ be a sufficiently large integer such that $N=n\mod2$. Now 
$u:\kk^n@>>>\kk^n$ has $i_e$ Jordan blocks of size $e$ ($e=1,2,3,\do$). Here $i_2,i_4,i_6,\do$ are even. Let 
$\D=\{e\in\{1,3,5,\do\};i_e\ge1\}$. If $\D=\em$ then $A(u)=\{1\}$, $\bA(u)=\{1\}$ and $\cg_\cf=\{1\}$ so that the
result is trivial.

In the remainder of this subsection we assume that $\D\ne\em$. Then $A(u)$ can be identified in the standard way 
with the $F_2$-subspace $\cp_{ev}(\D)$ of $\cp(\D)$ and the group of characters $\hA(u)$ of $A(u)$ (which may be 
canonically identified with the $F_2$-vector space dual to $A(u)$) becomes $\bcp(\D)$; the obvious pairing 
$A(u)\T\hA(u)@>>>F_2$ is induced by the inner product $L,L'\m|L\cap L'|\mod2$ on $\cp(\D)$.

To the partition $1i_1+2i_2+3i_3+\do$ of $n$ we associate a pair $(A,B)$ as in \cite{\ICC, 11.7} (with $N,M$ 
replaced by $n,N$). We have $A=\{\oca_i;i\in[0,N]\cap2\NN\}$, $B=(\oca_i;i\in[0,N]\cap(2\NN+1)\}$ where 
$\oca_0\le\oca_1\le\oca_2\le\do\le\oca_N$ is obtained from a sequence $a_0\le a_1\le a_2\le\do\le a_N$ as in 2.1.
(Here we use that $C$ is special.) Now the definitions and results in \S2 are applicable. As in \cite{\ORA, 4.5} 
(if $N$ is even) or \cite{\ORA, 4.6} (if $N$ is odd) the family $\cf$ is in canonical bijection with $\fT'$ in 
2.6. 

We arrange the intervals in $\ci$ in increasing order $I_{(1)},I_{(2)},\do,I_{(f)}$ (any element in $I_{(1)}$ is 
smaller than any element in $I_{(2)}$, etc.). We arrange the elements of $\D$ in increasing order 
$e_1<e_2<\do<e_{f'}$; then $f=f'$ and we have a bijection $\ci\lra\D$, $I_{(h)}\lra e_h$; moreover we have 
$|I_{(h)}|=i_{e_h}$ for $h\in[1,f]$; see \cite{\ICC, 11.7}. Using this bijection we see that $A(u)$ is identified
with $\cp_{ev}(\ci)$ and $\hA(u)$ is identified with $\bcp(\ci)$. For any $X\in\cp_{ev}(\ci)$, the assignment 
$L\m|X\cap L|\mod2$ can be viewed as an element of $\bcp(\ci)_1^*$ (the dual space of $\bcp(\ci)_1$ in 2.7 which 
by 2.7(b) is an $F_2$-vector space of dimension $2^{\lfloor\mu/2\rfloor}$). This induces a (surjective) 
$F_2$-linear map $\p:\cp_{ev}(\ci)@>>>\bcp(\ci)_1^*$. We will show that 

(a) $\ker\p=\ck(u)$ ($\ck(u)$ as in 0.1).
\nl
We identify $\Irr_CW$ with $T'$ (see 2.7) via the restriction of the bijection in \cite{\ICC, (13.2.5)} if $N$ 
is odd or \cite{\ICC, (13.2.6)} if $N$ is even (we also use the description of the Springer correspondence in 
\cite{\ICC, 13.3}). Under this identification the subset $\Irr^*_CW$ of $\Irr_CW$ becomes the subset $T_1$ (see 
2.7) of $T'$. Via the identification $\bcp(\ci)'\lra T'$ in 2.7 and $\hA(u)\lra\bcp(\ci)$ (see above), the map
$E\m\cv_E$ from $T'$ to $\hA(u)$ becomes the obvious imbedding $\bcp(\ci)_0@>>>\bcp(\ci)$ (we use again 
\cite{\ICC, 13.3}). By definition, $\ck(u)$ is the set of all $X\in\cp_{ev}(\ci)$ such that for any 
$L\in\cp(\ci)$ representing a vector in $\bcp(\ci)_1$ we have $|X\cap L|=0\mod2$. Thus, (a) holds.

Using (a) we have canonically $\bA(u)=\bcp(\ci)_1^*$ via $\p$. We have an $F_2$-linear map 
$\bcp(\ci)_1@>>>\bcp(\cj)_0$ (see 2.6) induced by $I_s\m\{a_{i_s},a_{i_{s+1}}\}$ for $s\in[0,\mu-1]\cap2\NN$ 
($I_s$ as in 2.7). This is an isomorphism; it corresponds to the bijection 2.7(c) under the identification
$T_1\lra\bcp(\ci)_1$ in 2.7 and the identification $\fT_1\lra\bcp(\cj)_0$ in 2.6. Hence we can identify 
$\bcp(\ci)_1^*$ with $\bcp(\cj)_0^*$ and with $\bcp(\cj)_1$ (see 2.6(a)). We obtain an identification 
$\bA(u)=\bcp(\cj)_1$.

By \cite{\ORA, 4.6} we have $\XX_\cf=\bcp_{ev}(\cj)$. Using 2.6(a) we see that 
$\bcp(\cj)=M(\bcp(\cj)_1)=M(\bA(u))$ canonically so that 0.4 holds in our case. From the arguments above we see 
that in our case 0.5 follows from 2.7(c).

\subhead 3.4\endsubhead
In 3.5-3.9 we consider the case where $G$ is simple adjoint of exceptional type. In each case we list the 
elements of the set $\Irr_CW$ for each special unipotent class $C$ of $G$; the elements of $\Irr_CW-\Irr^*CW$ are
enclosed in $[\,]$. (The notation for the various $C$ is as in \cite{\SPAA}; the notation for the objects of
$\Irr W$ is as in \cite{\SPAA} (for type $E_n$) and as in \cite{\ORA, 4.10} for type $F_4$.) In each case the 
structure of $A(u),\bA(u)$ (for $u\in C$) is indicated; here $S_n$ denotes the symmetric group in $n$ letters.
The order in which we list the objects in $\Irr_CW$ corresponds to the following order of the irreducible 
representations of $A(u)=S_n$:

$1,\e$ ($n=2$); $1,r,\e$ ($n=3,G\ne G_2$); $1,r$ ($n=3,G=G_2$); $1,\l^1,\l^2,\s$  ($n=4$);

$1,\nu,\l^1,\nu',\l^2,\l^3$ ($n=5$)
\nl
(notation of \cite{\ORA, 4.3}). Now 0.4 and 0.5 follow in our case from the tables in 3.5-3.9 and the definitions
in \cite{\ORA, 4.8-4.13}. (In those tables $S_n$ is the symmetric group in $n$ letters.)

\subhead 3.5\endsubhead
Assume that $G$ is of type $E_8$.

$\Irr_{E_8}W=\{1_0\}$; $A(u)=\{1\},\bA(u)=\{1\}$

$\Irr_{E_8(a_1)}W=\{8_1\}$; $A(u)=\{1\},\bA(u)=\{1\}$

$\Irr_{E_8(a_2)}W=\{35_2\}$; $A(u)=\{1\},\bA(u)=\{1\}$

$\Irr_{E_7A_1}W=\{112_3,28_8\}$; $A(u)=S_2,\bA(u)=S_2$

$\Irr_{D_8}W=\{210_4,160_7\}$; $A(u)=S_2,\bA(u)=S_2$

$\Irr_{E_7(a_1)A_1}W=\{560_5,[50_8]\}$; $A(u)=S_2,\bA(u)=\{1\}$

$\Irr_{E_7(a_1)}W=\{567_6\}$; $A(u)=\{1\},\bA(u)=\{1\}$

$\Irr_{D_8(a_1)}W=\{700_6,300_8\}$; $A(u)=S_2,\bA(u)=S_2$

$\Irr_{E_7(a_2)A_1}W=\{1400_7,1008_9, 56_{19}\}$; $A(u)=S_3,\bA(u)=S_3$

$\Irr_{A_8}W=\{1400_8,1575_{10},350_{14}\}$; $A(u)=S_3,\bA(u)=S_3$

$\Irr_{D_7(a_1)}W=\{3240_9,[1050_{10}]\}$; $A(u)=S_2,\bA(u)=\{1\}$

$\Irr_{D_8(a_3)}W=\{2240_{10}, [175_{12}], 840_{13}\}$; $A(u)=S_3,\bA(u)=S_2$

$\Irr_{D_6A_1}W=\{2268_{10}, 1296_{13}\}$; $A(u)=S_2,\bA(u)=S_2$

$\Irr_{E_6(a_1)A_1}W=\{4096_{11},4096_{12}\}$; $A(u)=S_2,\bA(u)=S_2$

$\Irr_{E_6}W=\{525_{12}\}$; $A(u)=\{1\},\bA(u)=\{1\}$

$\Irr_{D_7(a_2)}W=\{4200_{12},3360_{13}\}$; $A(u)=S_2,\bA(u)=S_2$

$\Irr_{E_6(a_1)}W=\{2800_{13},2100_{16}\}$; $A(u)=S_2,\bA(u)=S_2$

$\Irr_{D_5A_2}W=\{ 4536_{13},[840_{14}]\}$; $A(u)=S_2,\bA(u)=\{1\}$

$\Irr_{D_6(a_1)A_1}W=\{6075_{14},[700_{16}]\}$; $A(u)=S_2,\bA(u)=\{1\}$

$\Irr_{A_6A_1}W=\{2835_{14}\}$; $A(u)=\{1\},\bA(u)=\{1\}$

$\Irr_{A_6}W=\{4200_{15}\}$; $A(u)=\{1\},\bA(u)=\{1\}$

$\Irr_{D_6(a_1)}W=\{5600_{15},2400_{17}\}$; $A(u)=S_2,\bA(u)=S_2$

$\Irr_{2A_4}W=\{4480_{16},4536_{18}, 5670_{18},1400_{20}, 1680_{22}, 70_{32}\}$; $A(u)=S_5,\bA(u)=S_5$

$\Irr_{D_5}W=\{2100_{20}\}$; $A(u)=\{1\},\bA(u)=\{1\}$

$\Irr_{(A_5A_1)''}W=\{5600_{21},2400_{23}\}$; $A(u)=S_2,\bA(u)=S_2$

$\Irr_{D_4A_2}W=\{4200_{15},[168_{24}]\}$; $A(u)=S_2,\bA(u)=\{1\}$

$\Irr_{A_4A_2A_1}W=\{2835_{22}\}$; $A(u)=\{1\},\bA(u)=\{1\}$

$\Irr_{A_4A_2}W=\{4536_{23}\}$; $A(u)=\{1\},\bA(u)=\{1\}$

$\Irr_{D_5(a_1)}W=\{2800_{25}, 2100_{28}\}$; $A(u)=S_2,\bA(u)=S_2$

$\Irr_{A_42A_1}W=\{4200_{24},3360_{25}\}$; $A(u)=S_2,\bA(u)=S_2$

$\Irr_{D_4}W=\{525_{36}\}$; $A(u)=\{1\},\bA(u)=\{1\}$

$\Irr_{A_4A_1}W=\{4096_{26}, 4096_{27}\}$; $A(u)=S_2,\bA(u)=S_2$

$\Irr_{A_4}W=\{2268_{30},1296_{33}\}$; $A(u)=S_2,\bA(u)=S_2$

$\Irr_{D_4(a_1)A_2}=\{2240_{28},840_{31}\}$; $A(u)=S_2,\bA(u)=S_2$

$\Irr_{A_3A_2}W=\{3240_{31},[972_{32}]\}$; $A(u)=S_2,\bA(u)=\{1\}$

$\Irr_{D_4(a_1)A_1}W=\{1400_{32},1575_{34}, 350_{38}\}$; $A(u)=S_3,\bA(u)=S_3$

$\Irr_{D_4(a_1)}W=\{1400_{37},1008_{39},56_{49}\}$; $A(u)=S_3,\bA(u)=S_3$

$\Irr_{2A_2}W=\{700_{42},300_{44}\}$; $A(u)=S_2,\bA(u)=S_2$

$\Irr_{A_3}W=\{567_{46}\}$; $A(u)=\{1\},\bA(u)=\{1\}$

$\Irr_{A_22A_1}W=\{560_{47}\}$; $A(u)=\{1\},\bA(u)=\{1\}$

$\Irr_{A_2A_1}W=\{210_{52},160_{55}\}$; $A(u)=S_2,\bA(u)=S_2$

$\Irr_{A_2}W=\{112_{63}, 28_{68}\}$; $A(u)=S_2,\bA(u)=S_2$ 

$\Irr_{2A_1}W=\{35_{74}\}$; $A(u)=\{1\},\bA(u)=\{1\}$

$\Irr_{A_1}W=\{8_{91}\}$; $A(u)=\{1\},\bA(u)=\{1\}$

$\Irr_{\em}W=\{1_{120}\}$; $A(u)=\{1\},\bA(u)=\{1\}$

\subhead 3.6\endsubhead
Assume that $G$ is adjoint of type $E_7$.

$\Irr_{E_7}W=\{1_0\}$; $A(u)=\{1\},\bA(u)=\{1\}$

$\Irr_{E_7(a_1)}W=\{7_1\}$; $A(u)=\{1\},\bA(u)=\{1\}$

$\Irr_{E_7(a_2)}W=\{27_2\}$; $A(u)=\{1\},\bA(u)=\{1\}$

$\Irr_{D_6A_1}W=\{56_3, 21_6\}$; $A(u)=S_2,\bA(u)=S_2$ 

$\Irr_{E_6}W=\{21_3\}$; $A(u)=\{1\},\bA(u)=\{1\}$

$\Irr_{E_6(a_1)}W=\{120_4,105_5\}$; $A(u)=S_2,\bA(u)=S_2$ 

$\Irr_{D_6(a_1)A_1}W=\{189_5,[15_7]\}$; $A(u)=S_2,\bA(u)=\{1\}$ 

$\Irr_{D_6(a_1)}W=\{210_6\}$; $A(u)=\{1\},\bA(u)=\{1\}$     

$\Irr_{A_6}W=\{105_6\}$; $A(u)=\{1\},\bA(u)=\{1\}$

$\Irr_{D_5A_1}W=\{168_6\}$; $A(u)=\{1\},\bA(u)=\{1\}$

$\Irr_{D_5}W=\{189_7\}$; $A(u)=\{1\},\bA(u)=\{1\}$

$\Irr_{D_6(a_2)A_1}W=\{315_7,280_9,35_{13}\}$; $A(u)=S_3,\bA(u)=S_3$  

$\Irr_{(A_5A_1)'}=\{405_8,189_{10}\}$; $A(u)=S_2,\bA(u)=S_2$ 

$\Irr_{D_5(a_1)A_1}W=\{378_9\}$; $A(u)=\{1\},\bA(u)=\{1\}$

$\Irr_{A_4A_2}W=\{210_{10}\}$; $A(u)=\{1\},\bA(u)=\{1\}$

$\Irr_{D_5(a_1)}W=\{420_{10},336_{11}\}$; $A(u)=S_2,\bA(u)=S_2$ 

$\Irr_{A_5''}W=\{105_{12}\}$; $A(u)=\{1\},\bA(u)=\{1\}$

$\Irr_{A_4A_1}W=\{512_{11},512_{12}\}$; $A(u)=S_2,\bA(u)=S_2$ 

$\Irr_{D_4}W=\{105_{15}\}$; $A(u)=\{1\},\bA(u)=\{1\}$

$\Irr_{A_4}W=\{420_{13},336_{14}\}$; $A(u)=S_2,\bA(u)=S_2$ 

$\Irr_{A_3A_2A_1}W=\{210_{13}\}$; $A(u)=\{1\},\bA(u)=\{1\}$

$\Irr_{A_3A_2}W=\{378_{14},[84_{15}]\}$; $A(u)=S_2,\bA(u)=\{1\}$ 

$\Irr_{D_4(a_1)A_1}W=\{405_{15},189_{17}\}$; $A(u)=S_2,\bA(u)=S_2$ 

$\Irr_{D_4(a_1)}W=\{315_{16},280_{18},35_{22}\}$; $A(u)=S_3,\bA(u)=S_3$ 

$\Irr_{(A_3A_1)''}W=\{189_{20}\}$; $A(u)=\{1\},\bA(u)=\{1\}$

$\Irr_{2A_2}W=\{168_{21}\}$; $A(u)=\{1\},\bA(u)=\{1\}$

$\Irr_{A_23A_1}W=\{105_{21}\}$; $A(u)=\{1\},\bA(u)=\{1\}$

$\Irr_{A_3}W=\{210_{21}\}$; $A(u)=\{1\},\bA(u)=\{1\}$

$\Irr_{A_22A_1}W=\{189_{22}\}$; $A(u)=\{1\},\bA(u)=\{1\}$

$\Irr_{A_2A_1}W=\{120_{25},105_{26}\}$; $A(u)=S_2,\bA(u)=S_2$ 

$\Irr_{3A_1''}W=\{21_{36}\}$; $A(u)=\{1\},\bA(u)=\{1\}$

$\Irr_{A_2}W=\{56_{30},21_{33}\}$; $A(u)=S_2,\bA(u)=S_2$ 

$\Irr_{2A_1}W=\{27_{37}\}$; $A(u)=\{1\},\bA(u)=\{1\}$

$\Irr_{A_1}W=\{7_{46}\}$; $A(u)=\{1\},\bA(u)=\{1\}$

$\Irr_{\em}W=\{1_{63}\}$; $A(u)=\{1\},\bA(u)=\{1\}$

\subhead 3.7\endsubhead
Assume that $G$ is adjoint of type $E_6$.

$\Irr_{E_6}W=\{1_0\}$; $A(u)=\{1\},\bA(u)=\{1\}$

$\Irr_{E_6(a_1)}W=\{6_1\}$; $A(u)=\{1\},\bA(u)=\{1\}$

$\Irr_{D_5}W=\{20_2\}$; $A(u)=\{1\},\bA(u)=\{1\}$

$\Irr_{A_5A_1}W=\{30_3,15_5\}$; $A(u)=S_2,\bA(u)=S_2$ 

$\Irr_{D_5(a_1)}W=\{64_4\}$; $A(u)=\{1\},\bA(u)=\{1\}$                    

$\Irr_{A_4A_1}W=\{60_5\}$; $A(u)=\{1\},\bA(u)=\{1\}$

$\Irr_{A_4}W=\{81_6\}$; $A(u)=\{1\},\bA(u)=\{1\}$

$\Irr_{D_4}W=\{24_6\}$; $A(u)=\{1\},\bA(u)=\{1\}$

$\Irr_{D_4(a_1)}W=\{80_7,90_8,20_{10}\}$; $A(u)=S_3,\bA(u)=S_3$ 

$\Irr_{2A_2}W=\{24_{12}\}$; $A(u)=\{1\},\bA(u)=\{1\}$

$\Irr_{A_3}W=\{81_{10}\}$; $A(u)=\{1\},\bA(u)=\{1\}$

$\Irr_{A_22A_1}W=\{60_{11}\}$; $A(u)=\{1\},\bA(u)=\{1\}$

$\Irr_{A_2A_1}w=\{64_{13}\}$; $A(u)=\{1\},\bA(u)=\{1\}$

$\Irr_{A_2}W=\{30_{15}, 15_{17}\}$; $A(u)=S_2,\bA(u)=S_2$ 

$\Irr_{2A_1}W=\{20_{20}\}$; $A(u)=\{1\},\bA(u)=\{1\}$

$\Irr_{A_1}W=\{6_{25}\}$; $A(u)=\{1\},\bA(u)=\{1\}$

$\Irr_{\em}W=\{1_{36}\}$; $A(u)=\{1\},\bA(u)=\{1\}$

\subhead 3.8\endsubhead
Assume that $G$ is of type $F_4$.

$\Irr_{F_4}W=\{1_1\}$; $A(u)=\{1\},\bA(u)=\{1\}$

$\Irr_{F_4(a_1)}W=\{4_2,2_3\}$; $A(u)=S_2,\bA(u)=S_2$ 

$\Irr_{F_4(a_2)}W=\{9_1\}$; $A(u)=\{1\},\bA(u)=\{1\}$

$\Irr_{B_3}W=\{8_1\}$; $A(u)=\{1\},\bA(u)=\{1\}$

$\Irr_{C_3}W=\{8_3\}$; $A(u)=\{1\},\bA(u)=\{1\}$

$\Irr_{F_4(a_3)}W=\{12_1,9_3,6_2,1_3\}$; $A(u)=S_4,\bA(u)=S_4$ 

$\Irr_{\tA_2}W=\{8_2\}$; $A(u)=\{1\},\bA(u)=\{1\}$

$\Irr_{A_2}W=\{8_4,[1_2]\}$; $A(u)=S_2,\bA(u)=\{1\}$ 

$\Irr_{A_1\tA_1}W=\{9_4\}$; $A(u)=\{1\},\bA(u)=\{1\}$

$\Irr_{\tA_1}W=\{4_5,2_2\}$; $A(u)=S_2,\bA(u)=S_2$ 

$\Irr_{\em}W=\{1_4\}$; $A(u)=\{1\},\bA(u)=\{1\}$

\subhead 3.9\endsubhead
Assume that $G$ is of type $G_2$.

$\Irr_{G_2}W$ is the unit representation; $A(u)=\{1\},\bA(u)=\{1\}$

$\Irr_{G_2(a_1)}W$ consists of the reflection representation and the one dimensional representation on which the 
reflection with respect to a long (resp.short) simple coroot acts nontrivially (resp. trivially); 
$A(u)=S_3,\bA(u)=S_3$ 

$\Irr_{\em}W=\{\sg\}$; $A(u)=\{1\},\bA(u)=\{1\}$

\subhead 3.10\endsubhead
This completes the proof of Theorem 0.4 and that of Corollary 0.5.

We note that the definition of $\cg_\cf$ given in \cite{\ORA} (for type $C_n,B_n$) is $\bcp(\cj)_1$ (in the 
setup of 3.2) and $\bcp(\cj)_0$ (in the setup of 3.3) which is noncanonically isomorphic to $\bA(u)$, unlike the
definition adopted here that is, $\bcp(\cj)_0$ (in the setup of 3.2) and $\bcp(\cj)_1$ (in the setup of 3.3) 
which makes $\cg_\cf$ canonically isomorphic to $\bA(u)$.

\head 4. Character sheaves \endhead
\subhead 4.1\endsubhead 
Let $\hG$ be a set of representatives for the isomorphism classes of character sheaves on $G$. For any conjugacy 
class $D$ in $G$ let $D_\o:=\{g_\o;g\in D\}$, a unipotent class in $G$. For any unipotent class $C$ in $G$ let 
$\cs_C$ be the set of conjugacy classes $D$ of $G$ such that $D_\o=C$. It is likely that the following property 
holds.

(a) {\it Let $K\in\hG$. There exists a unique unipotent class $C$ of $G$  such that 

-for any $D\in\cs_C$, $K|_D$ is a local system (up to shift);

-for some $D\in\cs_C$, we have $K|_D\ne0$;

-for any unipotent class $C'$ of $G$ such that $\dim C'\ge\dim C$, $C'\ne C$ and any $D\in\cs_{C'}$ we have 
$K|_D=0$. 
\nl
We say that $C$ is the unipotent support of $K$.}
\nl
(The uniqueness part is obvious.) Note that \cite{\SUP, 10.7} provides some support (no pun intended) for (a).

We shall now try to make (a) more precise in the case where $K\in\hG^{un}$, the subset of $\hG$ consisting of 
unipotent character sheaves (that is $\hG^{un}=\hG_{\bbq}$ with the notation of \cite{\INT, 4.2}). As in 
\cite{\INT, 4.6} we have a partition $\hG^{un}=\sqc_{\cf}\hG^{un}_{\cf}$ where $\cf$ runs over the families of 
$W$.

In the remainder of this section we fix a family $\cf$ of $W$ and we denote by $C$ the special unipotent class of
$G$ such that $E_C\in\cf$, see 0.1; let $u\in C$. Let $\G=\bA(u)$ and let $Z(u)@>j'>>A(u)@>h>>\G$ be the obvious 
(surjective) homomorphisms; let $j=hj':Z(u)@>>>\G$. Let $[\G]$ be the set of conjugacy classes in $A(u)$. For 
$D\in\cs_C$ let $\ph(D)$ be the conjugacy class of $j(g_s)$ in $\G$ where $g\in D$ is such that $g_\o=u$; 
clearly such $g$ exists and is unique up to $Z(u)$-conjugacy so that the conjugacy class of $j(g_s)$ is 
independent of the choice of $g$.  Thus we get a (surjective) map $\ph:\cs_C@>>>[\G]$. For $\g\in[\G]$ we set 
$\cs_{C,\g}=\ph\i(\g)$. We now select for each $\g\in[\G]$ an element $x_\g\in\g$ and we denote by 
$\Irr Z_\G(x_\g)$ a set of representatives for the isomorphism classes of irreducible representations of 
$Z_\G(x_\g):=\{y\in\G;yx_\g=x_\g y\}$ (over $\bbq$). Let $D\in\cs_{C,\g}$, $\ce\in\Irr Z_\G(x_\g)$. We can find 
$g\in D$ such that $g_\o=u,j(g_s)=x_\g$ (and another choice for such $g$ must be of the form $bgb\i$ where 
$b\in Z(u)$, $j(b)\in Z_\G(x_\g)$). Let $\ce^D$ be the $G$-equivariant local system on $D$ whose stalk at 
$g_1\in D$ is $\{z\in G;zgz\i=g_1\}\T\ce$ modulo the equivalence relation $(z,e)\si(zh\i,j(h)e)$ for all 
$h\in Z(g)$. If $g$ is changed to $g_1=bgb\i$ ($b$ as above) then $\ce^D$ is changed to the $G$-equivariant local
system $\ce^D_1$ on $D$ whose stalk at $g'\in D$ is $\{z'\in G;z'g_1z'{}\i=g'\}\T\ce$ modulo the equivalence 
relation $(z',e')\si(z'h'{}\i,j(h')e)$ for all $h'\in Z(g_1)$. We have an isomorphism of local systems 
$\ce^D@>\si>>\ce^D_1$ which for any $g'\in D$ maps the stalk of $\ce^D$ at $g'$ to the stalk of $\ce^D_1$ at $g'$
by the rule $(z,e)\m(zb\i,j(b)e)$. (We have $zb\i g_1bz\i=zgz\i=g'$.) This is compatible with the equivalence 
relations. Thus the isomorphism class of the local system $\ce^D$ does not depend on the choice of $g$.

The properties (b),(c) below appear to be true ($[\,]$ denotes a shift). 

(b) {\it Let $K\in\hG^{un}_{\cf}$. There exists a unique $\g\in[\G]$ and a unique $\ce\in\Irr Z_\G(x_\g)$ such 
that

(i) if $D\in\cs_{C,\g}$, we have $K|_D\cong\ce^D[\,]$;

(ii) if $D\in\cs_{C,\g'}$ with $\g'\in[\G]-\{\g\}$, we have $K|_D=0$;

(iii) for any unipotent class $C'$ of $G$ such that $\dim C'\ge\dim C$, $C'\ne C$ and any $D\in\cs_{C'}$ we have 
$K|_D=0$.}
\nl
(c) {\it $K\m(\g,\ce)$ in (b) defines a bijection $\hG^{un}_{\cf}@>\si>>M(\G)$.}
\nl
Note that (b)(iii) follows from \cite{\SUP, 10.7}, at least if $p$ is sufficiently large or $0$.

In the case where $G$ is of type $E_8$ and $\cf$ contains the irreducible representation of degree $4480$ (so that
$\G=S_5$), (b)(i),(b)(ii),(c) have been already stated (without proof) in \cite{\INT, 4.7}.

For any finite dimensional representation $E$ of $W$ (over $\bbq$) let $\uE$ be the intersection cohomology 
complex on $G$ with coefficients in the local system with monodromy given by the $W$-module $E$ on the open set 
of regular semisimple elements. We have an imbedding $\cf@>>>\hG^{un}_{\cf}$, $E\m\uE[\,]$. Composing this 
imbedding with the map $\hG^{un}_{\cf}@>\si>>M(\G)$ in (c) (which we assume to hold) we obtain an imbedding 
$\cf@>>>M(\G)$. We expect that:

(d) {\it The imbedding $\cf@>>>M(\G)$ defined above coincides with the imbedding $\cf@>>>M(\G)$ in 
\cite{\ORA, Sec.4}.}
\nl
Note that 0.6 can be regarded as evidence for the validity of (b),(c),(d). Further evidence is given in 4.2-4.5.

\subhead 4.2\endsubhead
Assume that $G$ is simply connected. Let $D\in\cs_C$. Let $s$ be a semisimple element of $G$ such that $su\in D$.
Let $C_0$ be the conjugacy class of $u$ in $Z(s)$. Let $W'$ be the Weyl group of $Z(s)$ regarded as a subgroup 
of $W$. For any finite dimensional $W'$-module $E'$ over $\bbq$ let $\uuE'$ be the intersection cohomology 
complex on $Z(s)$ defined in terms of $Z(s),E'$ in the same way as $\uE$ was defined in terms of $G,E$. Using 
\cite{\CSII, (8.8.4)} and the $W$-equivariance of the isomorphism in {\it loc.cit.} we see that:

(a) {\it $\uE|_{sC_0}\cong(\un{\un{E|_{W'}}})|_{sC_0}[\,]$.}
\nl
Now, if $K\in\hG^{un}_{\cf}$ is of the form $\uE[\,]$ for some $E\in\cf$ then the computation of $K|_D$ is 
reduced by (a) to the computation of $\uuE'|_{sC_0}$ for any irreducible $W'$-module $E'$ such that
$(E':E_{W'})>0$ (here $(E':E_{W'})$ is the multiplicity of $E'$ in $E|_{W'}$). If for such $E'$ we define a 
unipotent class $\cc_{E'}$ of $Z(s)$ by $E'\in\Irr_{\cc_{E'}}W'$ then, by a known property of $\uuE'$, we have 
(with notation of 0.1 with $G$ replaced by $Z(s)$): 

(b) if $C_0=\cc_{E'}$ then $\uuE'|_{sC_0}[\,]$ is the irreducible $Z(s)$-equivariant local system corresponding 
to $\cv_{E'}$; 

(c) if $C_0\ne\cc_{E'}$ and $\dim C_0\ge\dim\cc_{E'}$ then $\uuE'|_{sC_0}=0$.
\nl
We say that $D$ is $E$-negligible if for any $E'\in\Irr W'$ such that $(E':E|_{W'})>0$ we have 
$\dim C_0>\dim\cc_{E'}$. 

(d) We say that $D$ is $E$-relevant if 

-there is a unique $E'_0\in\Irr W'$ such that $(E'_0:E|_{W'})=1$ and $\cc_{E'}=C_0$ (we then write $E_!=E'_0$);

-for any $E'\in\Irr W'$ such that $(E':E|_{W'})>0,E'\ne E_!$ we have $\dim C_0>\dim\cc_{E'}$.
\nl
It is likely that $D$ is always $E$-negligible or $E$-relevant. If $D$ is 
$E$-negligible then $\uE|_{sC_0}=0$ (hence $K|_D=0$); if $D$ is $E$-relevant then $\uE|_{sC_0}$ (hence $K|_D$) 
can be explicitly computed using (b),(c). 

In the remainder of this subsection we assume in addition that $G$ is almost simple of exceptional type and that
$C$ is a distinguished unipotent class. In these cases one can verify that $D$ is $E$-negligible or $E$-relevant 
for any $E\in\cf$ hence $K|_D$ can be explicitly computed and we can check that 4.1(b) holds. Moreover, we can 
compute $K|_D$ for any $K\in\hG^{un}_{\cf}$ (not necessarily of form $\uE[\,]$) using an appropriate analogue of 
(a) (coming again from \cite{\CSII, (8.8.4)}) and the appropriate analogues of (b),(c) (given in \cite{\ICC}). We
see that 4.1(b) holds again. Moreover we see that 4.1(c),(d) hold in these cases.

\subhead 4.3\endsubhead
In this subsection we assume that $G$ is of type $E_8$ and $C$ is distinguished. In this subsection we indicate
for each $D\in\cs_C$ the set $\cf_D=\{E\in\cf; D\text{ is }E-\text{\rm relevant}\}$ and we describe the map 
$E\m E_!$ (see 
4.2(d)). (Note that if $E\in\cf-\cf_D$, $D$ is $E$-negligible.) The notation is as in \cite{\SPAA}. 
We denote by $g_i$ an element of order $i$ of $A(u)$ (except that if $A(u)=S_5$, $g_2$ denotes a 
transposition and we denote by $g'_2$ an element of $A(u)$ whose centralizer has order $8$). For each $g_i$ we
denote by $\dg_i$ a semisimple element of $Z(u)$ that represents $g_i$; similarly when $A(u)=S_5$, we denote by 
$\dg'_2$ a semisimple element of $Z(u)$ that represents $g'_2$.
We write $\cf_{g_i}$ (resp. $\cf_{g'_2}$) instead of $\cf_D$ where $D$ is the $G$-conjugacy class of $u\dg_i$ 
(resp. of $u\dg'_2$). We write $\ch_{g_i}$ (resp. $\ch_{g'_2}$) for the set of all $E_!\in\Irr W'$ where $E$ runs 
through $\cf_{g_i}$ (resp. $\cf_{g'_2}$); here $W'\sub W$ is the Weyl group of $Z(\dg_i)$ (resp. $Z(\dg'_2)$) and 
$E_!$ is as in 4.2(d). We write $C_{g_i}$ (resp. $C_{g'_2}$) for the conjugacy class of $u$ in $Z(\dg_i)$ (resp.
$Z(\dg'_2)$). 

Assume that $C$ is the regular unipotent class. Then $A(u)=\{1\}$, $\cf_{g_1}=\ch_{g_1}=\{1_0\}$.

Assume that $C$ is the subregular unipotent class. Then $A(u)=\{1\}$, $\cf_{g_1}=\ch_{g_1}=\{8_1\}$.

Assume that $C=E_8(a_2)$. Then $A(u)=\{1\}$, $\cf_{g_1}=\ch_{g_1}=\{35_2\}$.

Assume that $C=E_7A_1$. Then $A(u)=S_2$, $Z(\dg_1)=G$, $Z(\dg_2)$ is of type $E_7A_1$,
$\cf_{g_1}=\ch_{g_1}=\{112_3,28_8\}$, $\cf_{g_2}=\{84_4\}$, $\ch_{g_2}=\{1_0\}$, $C_{g_2}=$ regular unipotent 
class.

Assume that $C=D_8$. Then $A(u)=S_2$, $Z(\dg_1)=G$, $Z(\dg_2)$ is of type $D_8$,
$\cf_{g_1}=\ch_{g_1}=\{210_4,160_7\}$, 
$\cf_{g_2}=\{50_8\}$, $\ch_{g_2}=\{1\}$, $C_{g_2}=$ regular unipotent class.

Assume that $C=E_7(a_1)A_1$. Then $A(u)=S_2$, $Z(\dg_1)=G$, $Z(\dg_2)$ is of type $E_7A_1$,
$\cf_{g_1}=\ch_{g_1}=\{560_5\}$,  $\cf_{g_2}=\{560_5\}$,  $\ch_{g_2}=\{7_1\bxt 1\}$,  
$C_{g_2}=$ subregular unipotent class in $E_7$ factor times regular unipotent class in $A_1$ factor.

Assume that $C=D_8(a_1)$. Then $A(u)=S_2$, $Z(\dg_1)=G$, $Z(\dg_2)$ is of type $D_8$,
$\cf_{g_1}=\ch_{g_1}=\{700_6,300_8\}$, $\cf_{g_2}=\{400_7\}$, $\ch_{g_2}=\{\text{\rm reflection repres.}\}$, 
$C_{g_2}=$ subregular unipotent class.

Assume that $C=E_7(a_2)A_1$. Then $A(u)=S_3$, $Z(\dg_1)=G$, $Z(\dg_2)$ is of type $E_7A_1$,
$Z(\dg_3)$ is of type $E_6A_2$, $\cf_{g_1}=\ch_{g_1}=\{1400_7,1008_9,56_{19}\}$, 
$\cf_{g_2}=\{1344_8\}$, $\ch_{g_2}=\{27_2\bxt1\}$, $C_{g_2}=$ subsubregular unipotent class in $E_7$-factor
times regular unipotent class in $A_1$ factor, 
$\cf_{g_3}=\{448_9\}$, $\ch_{g_3}=\{1\}$, $C_{g_3}=$ regular unipotent class.

Assume that $C=A_8$. Then $A(u)=S_3$, $Z(\dg_1)=G$, $Z(\dg_2)$ is of type $D_8$,
$Z(\dg_3)$ is of type $A_8$, $\cf_{g_1}=\ch_{g_1}=\{1400_8,1575_{10},350_{14}\}$, 
$\cf_{g_2}=\{1050_{10}\}$, $\ch_{g_2}=\{28-\text{dimensional repres}.\}$, $C_{g_2}=$ unipotent class with
Jordan blocks of size $5,11$, 
$\cf_{g_3}=\{175_{12}\}$, $\ch_{g_3}=\{1\}$, $C_{g_3}=$ regular unipotent class.

Assume that $C=D_8(a_3)$. Then $A(u)=S_3$, $Z(\dg_1)=G$, $Z(\dg_2)$ is of type $D_8$, $Z(\dg_3)$ is of type 
$E_6A_2$, $\cf_{g_1}=\ch_{g_1}=\{2240_{10},840_{13}\}$, 
$\cf_{g_2}=\{1400_{11}\}$, $\ch_{g_2}=\{56-\text{dimensional repres}.\}$, $C_{g_2}=$ unipotent class with
Jordan blocks of size $7,9$, $\cf_{g_3}=\{2240_{10}\}$, $\ch_{g_3}=\{6_1\bxt1\}$, 
$C_{g_3}=$ subregular unipotent class in $E_6$-factor times regular unipotent class in $A_1$ factor.

Assume that $C=2A_4$. Then $A(u)=S_5$, 

$Z(\dg_1)=G$, $Z(\dg_2)$ is of type $E_7A_1$, $Z(\dg'_2)$ is of type $D_8$, 

$Z(\dg_3)$ is of type $E_6A_2$, $Z(\dg_4)$ is of type $D_5A_3$, $Z(\dg_5)$ is of type $A_4A_4$, 

$Z(\dg_6)$ is of type $A_5A_2A_1$,

$\cf_{g_1}=\ch_{g_1}=\{4480_{16},4536_{18}, 5670_{18},1400_{20}, 1680_{22}, 70_{32}\}$,

$\cf_{g_2}=\{7168_{17},5600_{19},448_{25}\}$, $\ch_{g_2}=\{315_7\ot1, 280_9\ot1,35_{13}\ot1\}$,

$C_{g_2}=D_6(a_1)A_1$ in $E_7$-factor times regular unipotent class in $A_1$-factor,

$\cf_{g'_2}=\{4200_{18},2688_{20}\e'',168_{24}\}$, 

$\ch_{g'_2}=\{\text{ repres. with symbol }(2<5;0<3),(2<3;0<5),(0<1,4<5)\}$,

$C_{g'_2}=$unipotent class with Jordan blocks of sizes $1,3,5,7$,

$\cf_{g_3}=\{3150_{18},1134_{20}\}$, $\ch_{g_3}=\{30_3\bxt1,15_5\bxt1\}$,

$C_{g_3}=A_5A_1$ in $E_6$-factor times regular unipotent class in $A_2$-factor,

$\cf_{g_4}=\{1344_{19}\}$, $\ch_{g_4}=\{5-\text{dimensional repres.}\}$, 

$C_{g_4}=$subregular unipotent class in $D_5$-factor times regular unipotent class in $A_3$-factor,

$\cf_{g_5}=\{420_{20}\}$, $\ch_{g_5}=\{1\}$, $C_{g_5}=$regular unipotent class,

$\cf_{g_6}=\{2016_{19}\}$, $\ch_{g_6}=\{1\}$, $C_{g_6}=$regular unipotent class.

In each case the $i$-th member of a list $\cf_?$ and the $i$-th member of the corresponding list $\ch_?$ are
related by the map $E\m E_!$. Note that the members of the list $\cf_{g'_2}$ (when $C=2A_4$) are not all in the
same family. But in all cases, the members of the list $\cf_g$ form exactly the subset of $\Irr W'$ corresponding
to the unipotent class $C_g$ under Springer's correspondence for $Z(\dg)$; 
thus they can be indexed by certain irreducible
representations of the group of components of the centralizer of $u$ in $Z(\dg)$ modulo the centre of $Z(\dg)$.
(Here $g$ is $g_i$ or $g'_2$.) From this one recovers the imbedding $\cf@>>>M(\bA(u))$ in geometric terms.

\subhead 4.4\endsubhead
In this subsection we assume that $G=Sp_4(\kk)$ and that $\cf$ is the family in $\Irr W$ containing the
reflection representation so that $C$ is the subregular unipotent class in $G$. Let $D$ be the conjugacy class in 
$G$ containing $sv$ where $s$ is semisimple with $Z(s)\cong SL_2(\kk)\T SL_2(\kk)$ and $v$ is a regular unipotent
element of $Z(s)$ so that $v\in C$. Let $D'$ be a conjugacy class in $G$ containing $s'v'$ where $s'$ is 
semisimple with $Z(s')\cong GL_2(\kk)$ and $v'$ is a regular unipotent element of $Z(s')$ so that $v'\in C$. In 
this case $\hG^{un}_{\cf}$ consists of four character sheaves $K_1,K_2,K_3,K_4$, the last one being cuspidal.
They can be characterized as follows.

$K_1|_C=\bbq[\,]$, $K_1|_D=0$, $K_1|_{D'}=\bbq[\,]$;

$K_2|_C=\cl[\,]$, $K_2|_D=0$, $K_2|_{D'}=\bbq[\,]$;

$K_3|_C=0$, $K_3|_D=\bbq[\,]$, $K_3|_{D'}=0$;

$K_4|_C=0$, $K_4|_D=\cl'[\,]$, $K_4|_{D'}=0$.
\nl
Here $\cl$ is a notrivial $G$-equivariant local system of rank $1$ on $C$, $\cl'$ is the inverse image of $\cl$
under the obvious map $D@>>>C$. We see that 4.1(b) holds for all $K\in\hG^{un}_{\cf}$ and 4.1(c),(d) hold. 

\subhead 4.5\endsubhead
Assume that $\cf$ is the family containing the unit representation of $W$. Then $C$ is the regular unipotent class
of $G$ and $\hG^{un}_{\cf}$ consists of a single character sheaf, namely $\bbq[\,]$. Clearly, 4.1(b),(c),(d) hold 
in this case.

Next we assume that $\cf$ is the family containing the sign representation of $W$. Then $C=\{1\}$ and 
$\hG^{un}_{\cf}$ consists of a single character sheaf, namely $K=\un{\sg}[\,]$. Note that for any semisimple 
class $D$ of $G$ we have $K|_D=\bbq[\,]$ so that 4.1(b),(c),(d) hold in this case.

\enddocument